\documentclass[12pt,a4paper,final]{amsart}

\usepackage{amssymb,amsmath,amsthm}
\usepackage{color}
\usepackage{enumerate}
\usepackage{rotating}
\usepackage{fullpage}
\usepackage{comment}
\usepackage{hyperref}
\DeclareMathOperator{\inte}{int}

\newtheorem{theorem}{Theorem}
\newtheorem{corollary}[theorem]{Corollary}

\theoremstyle{definition}
\newtheorem{definition}{Definition}

\title{A computer-assisted proof of the existence of {S}male horseshoe for the folded-towel map}

\author{Anna Gierzkiewicz}
\email{\href{mailto:anna.gierzkiewicz@urk.edu.pl}{anna.gierzkiewicz@urk.edu.pl}}
\address{Department of Applied Mathematics, University of Agriculture in Krak\'ow,
ul. Balicka 253c,
30--198 Krak\'ow, Poland
}

\begin{document}

\begin{abstract}

The paper contains a rigorous proof of existence of symbolic dynamics chaos in the generalized H\'enon map's 4th iterate $H^4$, which was conjectured in the paper \textit{A 3{D} {S}male Horseshoe in a Hyperchaotic Discrete-Time System} of Li and Yang, 2007. We prove also the uniform hyperbolicity of the invariant set with symbolic dynamics. The proofs are computer-assisted with the use of C++ library \textit{CAPD} for interval arithmetic, differentiation and integration.
%\medskip

%\noindent \textbf{Keywords:} 
%\medskip

%\noindent \textbf{MSC Classification:} 
\end{abstract}

\maketitle

%----------------------------------------------------------------------------------------
%	Intro
%----------------------------------------------------------------------------------------

\section{Introduction}

The H\'enon map \cite{henon1976}
 is a well known example of chaotic diffeomorphism on $\mathbb{R}^2$, and has been widely developed and generalized to many different contexts. R\"ossler in 1979 \cite{R} introduced a~similar $\mathbb{R}^3$ map
to study hyper-chaos, which he described as `a higher form of chaos with two directions of hyperbolic instability on the attractor'. The system gained the name of `folded-towel map', thanks to its attractor's shape. Further generalization to $\mathbb{R}^n$ and its study via Lyapunov characteristic exponents may be found in \cite{Baier}.
For other studies of 3D H\'enon-like maps by estimating their two maximal Lyapunov exponents see \cite{3DHenonySimo}, also \cite{3DHenonyRichter}.

\begin{figure}[h]
\includegraphics[width=0.45\textwidth]{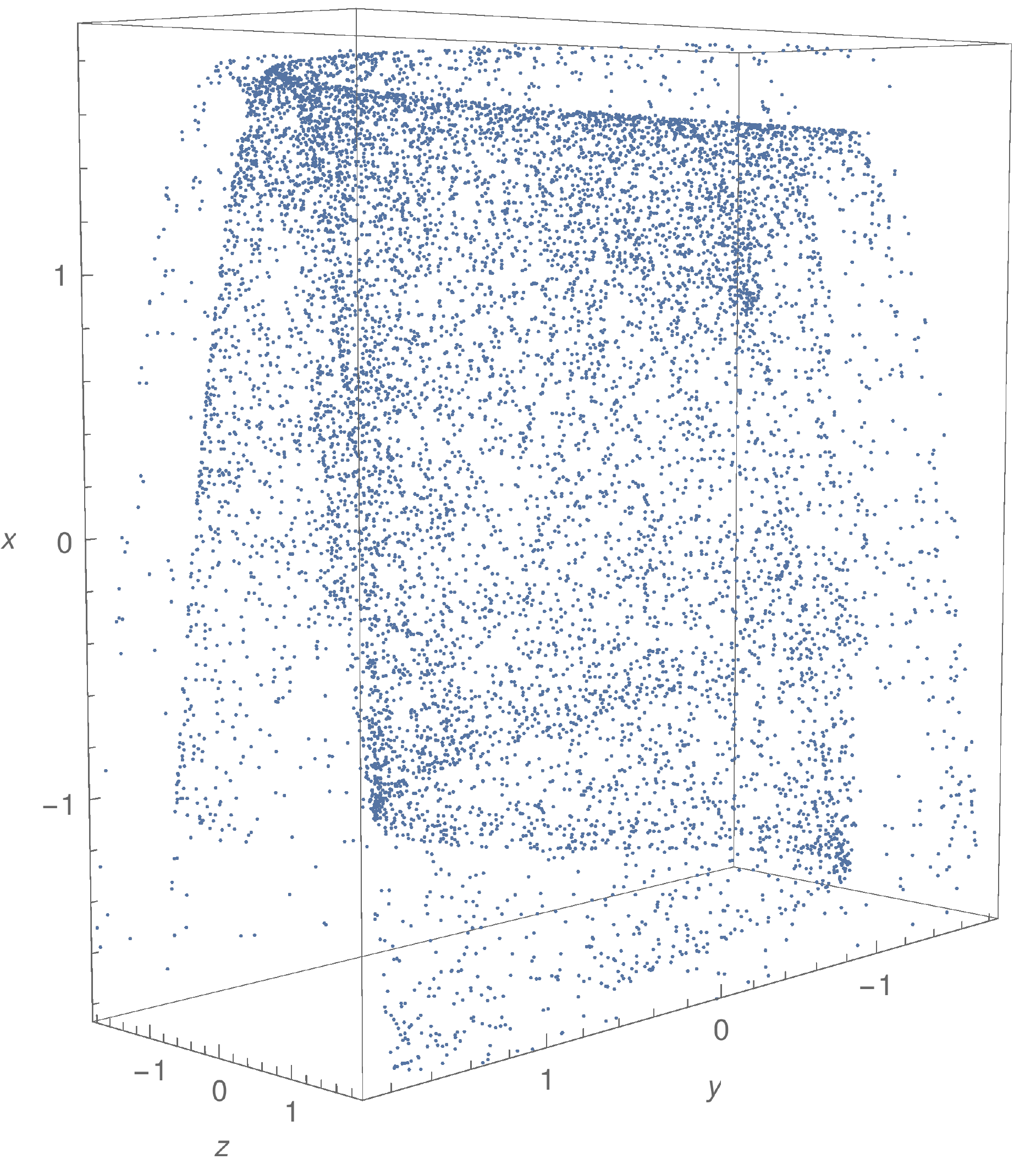}
\hfil 
\includegraphics[width=0.45\textwidth]{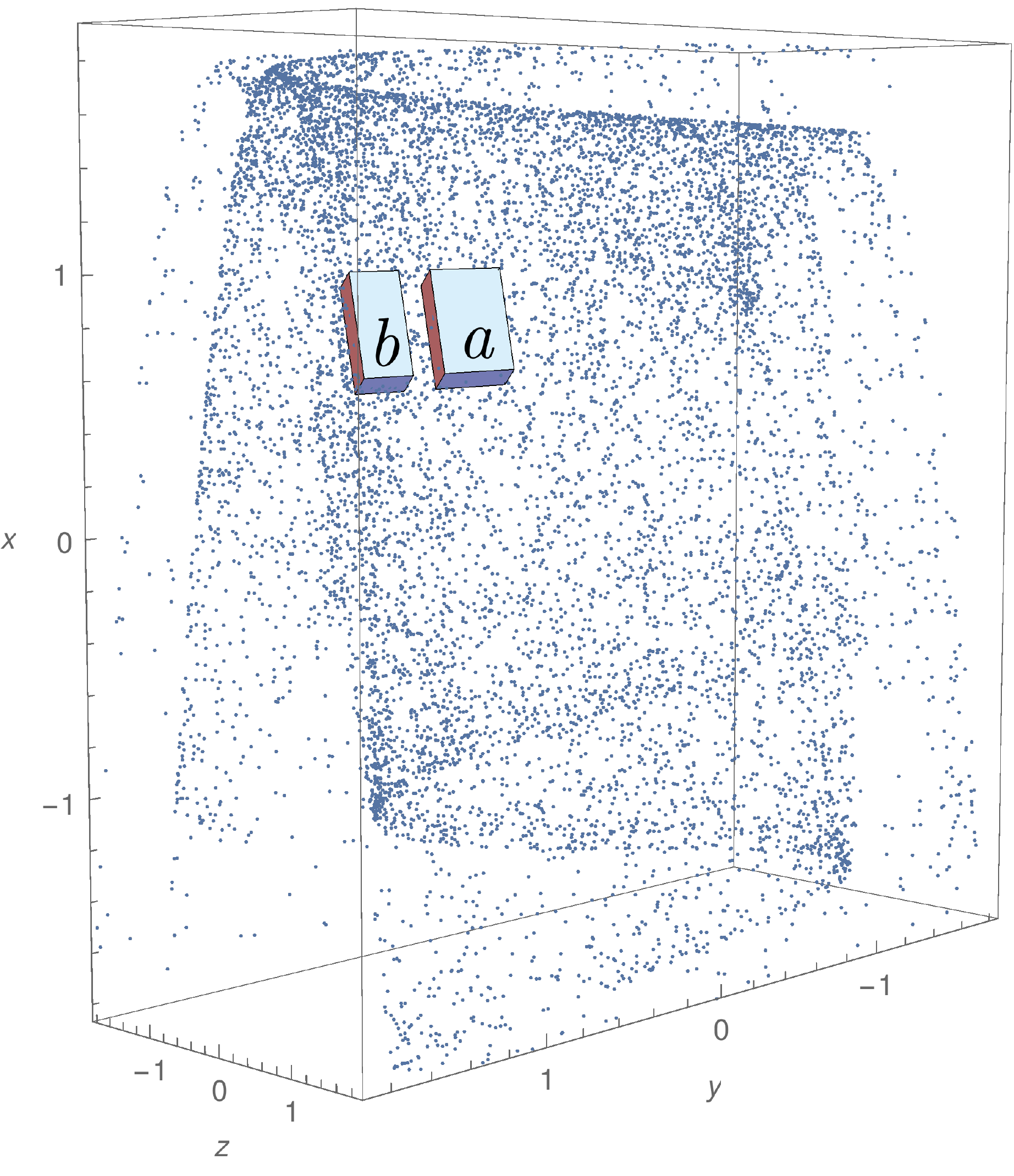}
\caption{\label{fig:towel}The `folded towel' attractor for H\'enon 3D map $H$ defined by Eq.\ \eqref{eq:h3d}.
\newline
To the right: The location of the sets $a$, $b$ on the attractor.}
\end{figure}

Our motivation is the article \cite{LiYang}, where 
the 3D case with fixed parameters is considered:
\begin{equation}\label{eq:h3d}
H(x,y,z) = (1.76 - y^2-0.1 z,\ x,\ y).
\end{equation} 
The system \eqref{eq:h3d} is also investigated, implemented as an electronic circuit, in \cite{Grassi}.
Its `folded towel' attractor is depicted on Fig. \ref{fig:towel}.
The authors of \cite{LiYang} observe that the hyper-chaos may be also studied as containing a~3D generalization of Smale horseshoe dynamics on a compact subset of $\mathbb{R}^3$. They show explicitly two cuboids $a$, $b$, which most probably contain a Smale horseshoe with two expanding directions for the fourth iterate of the map (see Fig. \ref{fig:towel}, to the right).

\begin{figure}[h]
\includegraphics[height=8.5cm]{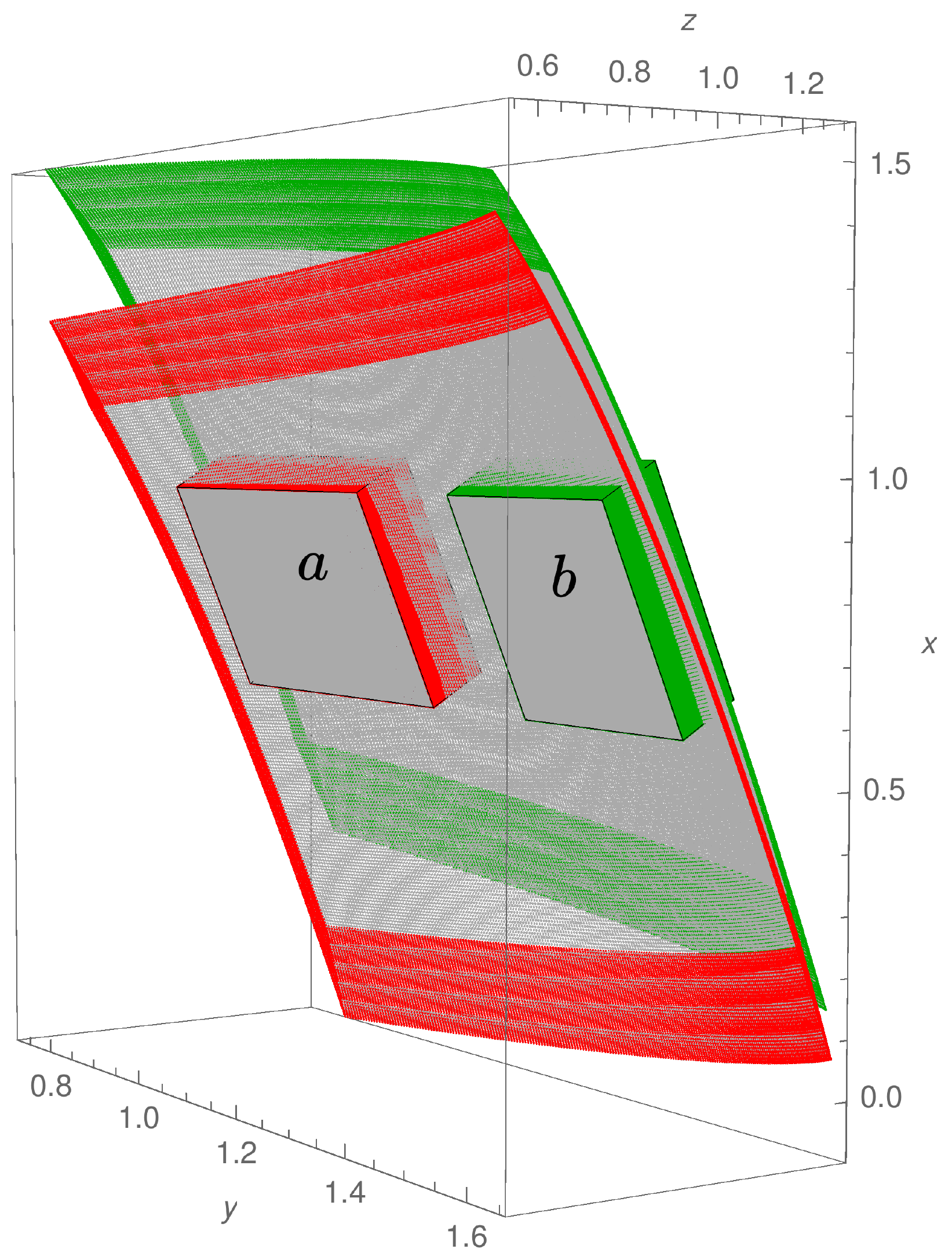}
\hfil
\includegraphics[height=8.5cm]{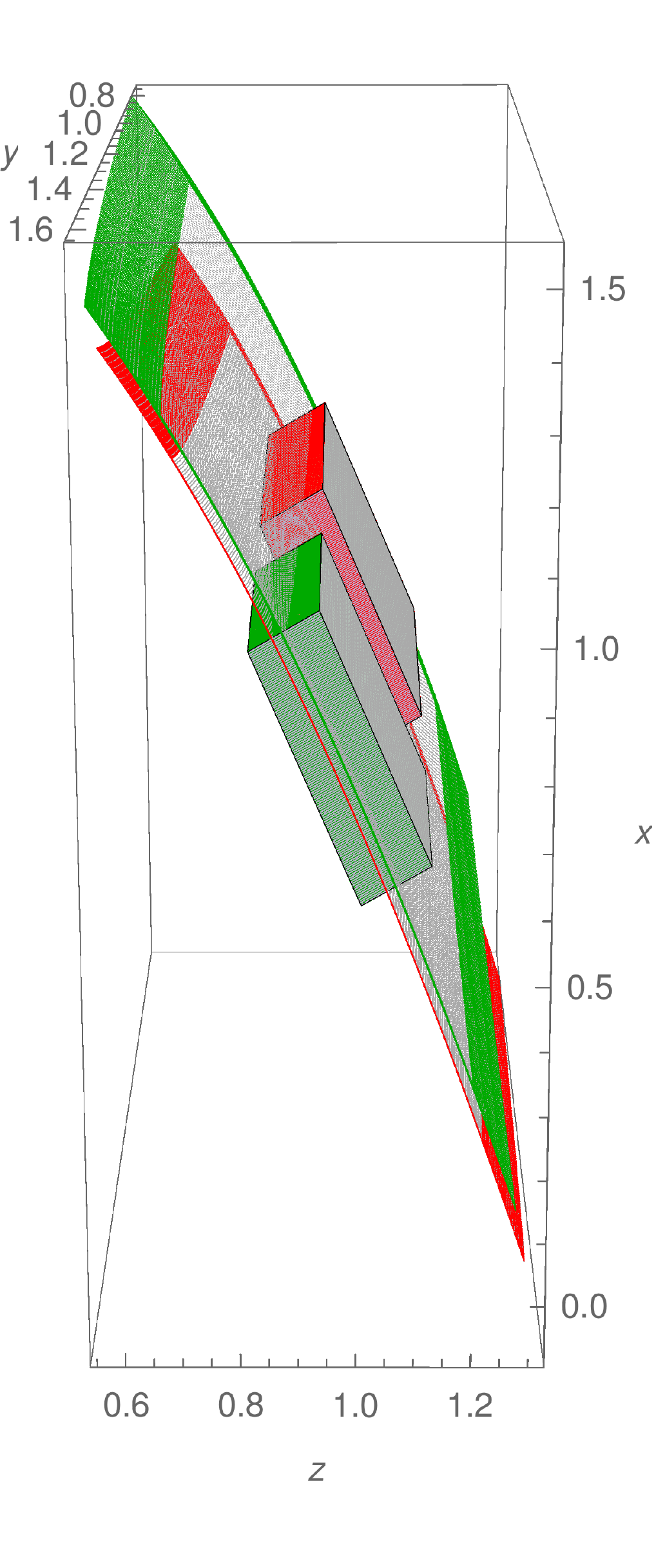}
\caption{\label{fig:H4ab}The location of the images $H^4(a)$, $H^4(b)$ with respect to sets $a$ and $b$ --- three-quarter view to the left and side view to the right. The exit sets and their images are marked in red (for $a$) and green (for $b$).}
\end{figure}

The images of $a$, $b$ via $H^4$ have properties that can be described intuitively as follows:
\begin{itemize}
	\item each image $H^4(a)$ and $H^4(b)$ intersects both $a$ and $b$;
	\item each image is compressed in the direction along the shortest edge of $a$ or $b$. This direction is `locally normal' to the attractor;	
	\item each image is expanded in two other directions, almost along the other edges of $a$ and $b$. 
\end{itemize}
Fig. \ref{fig:H4ab} shows the way the images $H^4(a)$ and $H^4(b)$ intersect $a$ and $b$. 

This is a well known method for establishing (hyper-)chaos by computer assisted proof, which can be implemented with interval arithmetic (see, for example, \cite{DW09, Barrio}). 
Our paper presents such a proof of the fact that the set $a\cup b$ is indeed a topological horseshoe, which implies the existence of symbolic dynamics for the system \eqref{eq:h3d} (Theorem \ref{th:symbolic}). Additionally we prove that the invariant set contained in $a\cup b$ is uniformly hyperbolic (Theorem \ref{th:myHyperbolic}).

The proofs of Theorems \ref{th:symbolic}, \ref{th:myHyperbolic} are computer-assisted. It means that their essential parts are C++ codes available on-line \cite{proof}. We are aware that such codes get outdated quickly and after a few years they may not compile with the current C++ compilers any more. We hope, however, that outlines of the proofs will help to understand the idea and the general construction of the codes.

\section{Topological covering and periodic orbits}\label{sec:tools}

For full description of h-sets and their covering relations see \cite{WZ1}. In \cite{AGPZ} one may also find a 2D simplified version, which may be useful for understanding the general idea. Here we introduce a short collection of necessary notions.

We use the standard notation for the closure, interior, and boundary of a topological set $A\subset\mathbb{R}^k$, which are  $\overline{A}$, $\inte A$, and $\partial A$, respectively. The $k$-dimensional open unit ball centred at the origin is denoted by $\mathbb{B}_k$. We use the balls in the maximum norm (cubes), for they are easily interpreted in the interval algebra as interval vectors:
\[
\mathbb{B}_k = \{\mathbf{x}\in \mathbb{R}^k \colon \|\mathbf{x}\|_{\infty}<1\}.
\]

\subsection{H-sets}
The basic object we work on is
\begin{definition}[\cite{WZ1}, Def. 3.1]
\emph{An h-set} is a quadruple $N=(|N|, u(N), s(N), C_N)$, where $|N|$ is a compact subset of $\mathbb{R}^n$, which we will call the \emph{support of an h-set} and
\begin{enumerate}
\item two numbers $u(N)$, $s(N) \in \mathbb{N}\cup\{0\}$ complement the dimension of space:
\[u(N) + s(N)=n\text{;}\]
we will call them the \emph{exit} and \emph{entry dimension}, or \emph{unstable} and \emph{stable dimension}, respectively;
\item the homeomorphism $C_N : \mathbb{R}^n\to\mathbb{R}^n=\mathbb{R}^{u(N)}\times \mathbb{R}^{s(N)}$ is such that
\[
C_N(|N|)=\overline{\mathbb{B}_{u(N)}}\times\overline{\mathbb{B}_{s(N)}}\text{.}
\]
\end{enumerate}
The support of an h-set is sometimes called simply an h-set, if it does not lead to confusion or its structure is unimportant. We also often use the notation `$f(N)$' interchangeably with `$f(|N|)$' to simplify formulas.
\end{definition}

Let us set also some useful notations:
\begin{align*}
\dim N &= n\text{,}\\
N_c &= \overline{\mathbb{B}_{u(N)}}\times\overline{\mathbb{B}_{s(N)}}\text{,}\\
N_c^- &= \partial\mathbb{B}_{u(N)}\times\overline{\mathbb{B}_{s(N)}}\text{,}\\
N_c^+ &= \overline{\mathbb{B}_{u(N)}}\times\partial\mathbb{B}_{s(N)}\text{,}\\
N^- &= C_N^{-1}(N_c^-)\text{,} \qquad N^+ = C_N^{-1}(N_c^+)\text{.}
\end{align*}

Therefore, we can assume that an h-set is a product of two unitary balls moved to some coordinate system with the exit set $N^-$ and entrance set $N^+$ distinguished.
The notions with the subscript $_c$ refer to the `straight' coordinate system in the image of $C_N$.

\subsection{Covering relation}
\begin{figure}[h]
	\includegraphics[height=5cm]{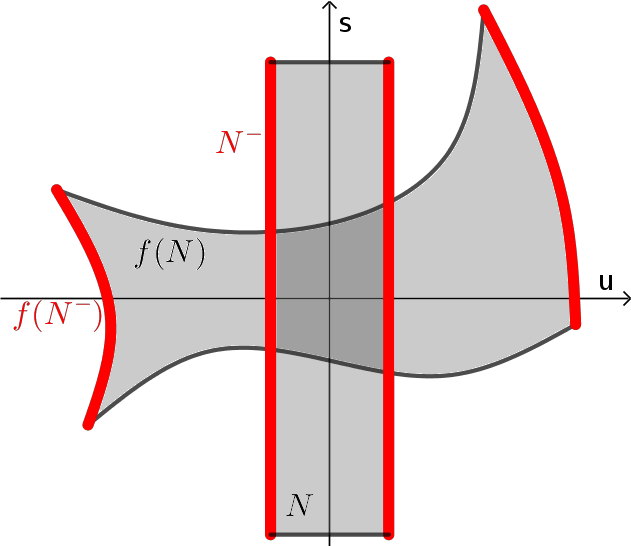}
	\hfil
	\includegraphics[height=5cm]{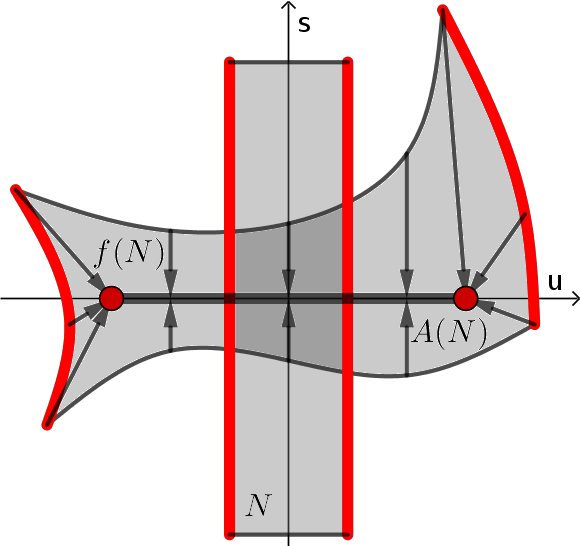}
	\caption{\label{fig:covering}To the left: an example of topological self-covering $N \stackrel{f}{\Longrightarrow}N$ in $\mathbb{R}^2$. To the right: $f(N)$ is homotopy equivalent to the image of $N$ via a linear map $A$ with certain properties (see Definition \ref{def:cov}). The exit sets and their images are marked in red.}
\end{figure}

\begin{figure}[h]
	\includegraphics[height=5cm]{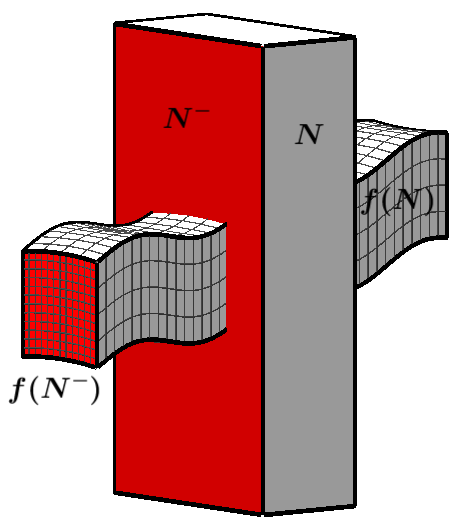}
	\hfil 	
	\includegraphics[height=5cm]{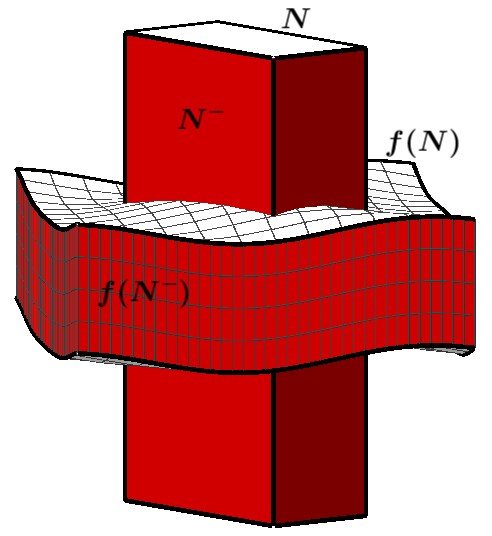}
	\caption{\label{fig:covering2}Examples of topological self-covering $N \stackrel{f}{\Longrightarrow}N$ in $\mathbb{R}^3$: with one exit direction (left) and two exit directions (right). The exit sets and their images are marked in red.}
\end{figure}

We define the topological covering:

\begin{definition}[\cite{WZ1}, Def. 3.4, simplified]\label{def:cov}
Let two h-sets $M$, $N$ be such that
$u(M)=u(N)=u$ and $s(M)=s(N) = s$. For a continuous map $f:|M|\to \mathbb{R}^n$ denote $f_c=C_N \circ f \circ C_M^{-1} : M_c \to \mathbb{R}^u\times \mathbb{R}^s$.

We say that $M$ \emph{$f$-covers} the h-set $N$ (denoted shortly  as  $M \stackrel{f}{\Longrightarrow}N$) if there exists a continuous homotopy $h:[0,1]\times M_c \to \mathbb{R}^u\times \mathbb{R}^s$, such that:
\begin{enumerate}
\item $h_0 = f_c $,
\item $h([0,1],M_c^-)\cap N_c = \varnothing$ \quad \textit{(the exit condition),}
\item $h([0,1],M_c)\cap N_c^+ = \varnothing$ \quad \textit{(the entry condition).}
\item Additionally, if $u>0$, then there exists a linear map $A:\mathbb{R}^u \to \mathbb{R}^u$ such that
\begin{align*}
h_1(x,y) &= (A(x),0) \quad \text{ for } x\in\overline{\mathbb{B}_u} \text{ and } y\in\overline{\mathbb{B}_s} \text{,}\\
A(\partial \mathbb{B}_u) &\subset \mathbb{R}^u\setminus\mathbb{B}_u.
\end{align*}
\end{enumerate}
\end{definition}
For some geometrical intuition of covering in low-dimensional cases, see Figs. \ref{fig:covering}, \ref{fig:covering2}.

Topological covering has a property of tracking orbits \cite{Mir}. In other words, for a chain covering $N_1 \stackrel{f}{\Longrightarrow}N_2\stackrel{f}{\Longrightarrow}N_3$, one can find a point in $N_1$ that is mapped into $N_2$ and then to $N_3$. Moreover, one can prove the existence of a periodic orbit related to a closed sequence of covering relations:

\begin{theorem}[\cite{WZ1}, Theorem 3.6, simplified]\label{th:periodic}
Suppose there exists a sequence of h-sets $N_0$, \ldots, $N_n=N_0$, such that
\[
N_0 \stackrel{f}{\Longrightarrow} N_1 \stackrel{f}{\Longrightarrow}\ldots \stackrel{f}{\Longrightarrow} N_n = N_0\text{,}
\]
then there exists a point $x\in \inte |N_0|$, such that $f^k(x) \in \inte|N_{k}|$ for $k=0,1,\ldots,n$ and $f^n(x)=x$.
\end{theorem}

In particular, if $N_0\stackrel{f}{\Longrightarrow} N_0$, then there exists a stationary point for the map $f$, contained in $N_0$.

\section{Detecting symbolic dynamics via covering relations}\label{sec:chaos}
%----------------------------------------------------------------------------------------
%	Chaos with covering
%----------------------------------------------------------------------------------------

We shall prove that the map \eqref{eq:h3d} is chaotic in the sense of symbolic dynamics for $H^4$. First let us recall this notion:

\subsection{Symbolic dynamics}
Let $\Sigma_2=\{0,1\}^{\mathbb{Z}}$ be the set of bi-infinite sequences of two symbols understood as a compact space with the metric
\[
\text{for }l=\{ l_n\}_{n\in\mathbb{Z}}\text{,\quad  }l'=\{ l'_n\}_{n\in\mathbb{Z}}\text{, }\qquad
\operatorname{dist}\left(l,l'\right) = \sum_{n=-\infty}^{+\infty} \frac{|l_n-l'_n|}{2^{|n|}}\text{,}
\]
which induces the product topology. The homeomorphism $\sigma:\Sigma_2\to \Sigma_2$, given by
$$
 (\sigma (l))_n =l_{n+1}\text{,}
$$
is called the \emph{shift map}. It has many interesting topological properties, in particular: the existence of dense orbits, or the density of periodic orbits' set in the whole space. 

In our study, by the chaotic behaviour of a discrete dynamical system induced by a homeomorphism $f:X\to X$ we understand the existence of a compact set $I\subset X$ invariant for $f$ (or sometimes its higher iterate) such that $f|_I$ is semi-conjugate to $\sigma$, {\it i.e.} there exists a continuous surjection $g:I\to \Sigma_2$ such that
\[
g\circ f|_I = \sigma \circ g.
\]

In other words, $f$ admits on $I$ at least as rich dynamics as $\sigma$ on $\Sigma_2$. It means, in particular, that the topological entropy of $f$ is greater or equal $\log {2}$.
The system $(\Sigma_2,\sigma)$ or any system semi-conjugate to it is often called in literature \emph{a symbolic dynamics system} \cite{Morse}.

Symbolic dynamics is sometimes used as one of the definitions of chaotic dynamics, because the discrete system $(\Sigma_2,\sigma)$ defined above evinces all the typical chaotic phenomena as transitivity, density of periodic orbits set or sensitivity to initial conditions. Also, there exists a periodic orbit of any prescribed period \cite{Morse}.

\subsection{Topological horseshoe}

\begin{figure}[h]
	\includegraphics[height=4cm]{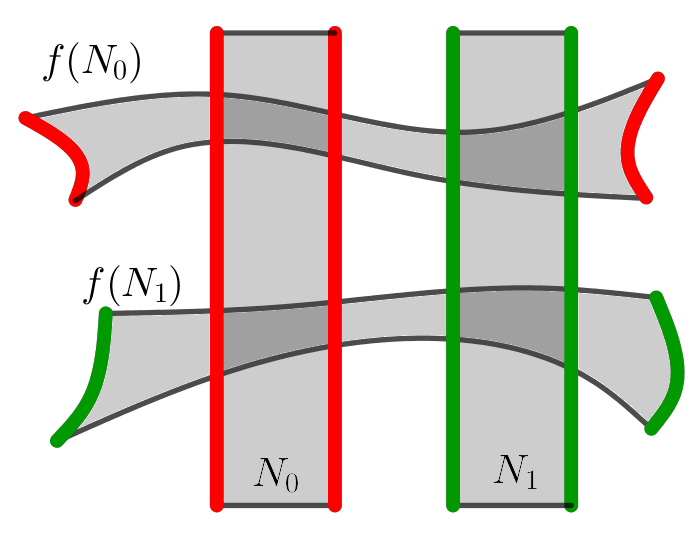}
	\hfil
	\includegraphics[height=4cm]{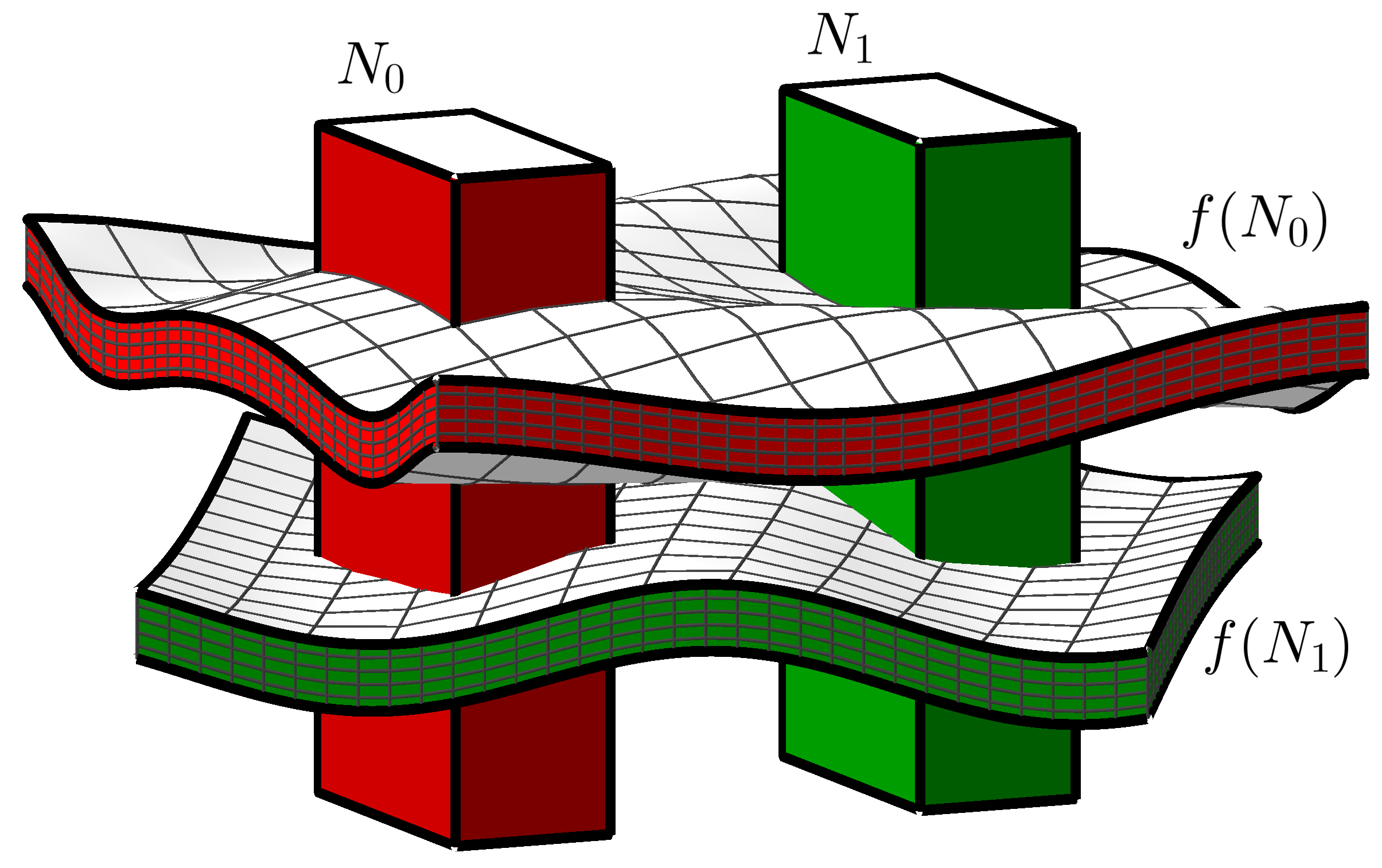}
	\caption{\label{fig:horseshoe}Topological horseshoes in $\mathbb{R}^2$ and $\mathbb{R}^3$: each $N_{0,1}$ covers itself and the other set. The exit sets of $N_0$ and $N_1$ are marked in red and green, respectively. Compare also to Fig. \ref{fig:H4ab}.}
\end{figure}

In general, to isolate a set with symbolic dynamics does not seem to be an easy task. It occurs, however, that some sets which fulfil certain covering relations must contain symbolic dynamics. An important example is:

\begin{definition}[Topological horseshoe]
Let $N_0$, $N_1 \subset \mathbb{R}^n$ be two disjoint h-sets. We say that a continuous map $f:\mathbb{R}^n \to \mathbb{R}^n$ is a \emph{topological horseshoe} for $N_0$, $N_1$ if (see Fig. \ref{fig:horseshoe})
\begin{equation}
\begin{array}{cc}
N_0 \overset{f}{\Longrightarrow}N_0\text{,} \quad
&N_0 \overset{f}{\Longrightarrow}N_1\text{,}
\\
N_1 \overset{f}{\Longrightarrow}N_0 \text{,} \quad
&N_1 \overset{f}{\Longrightarrow}N_1.
\end{array}
\end{equation}
\end{definition}

It can be shown that for any topological horseshoe we obtain symbolic dynamics.

\begin{theorem}[\cite{GZ}, Theorem 18]
Let $f$ be a topological horseshoe for $N_0$, $N_1$. Denote by $I = \operatorname{Inv}(N_0\cup N_1)$ the invariant part of the set $N_0\cup N_1$ under $f$, and define a map $g: I \to \Sigma_2$ by
\[
g(x)_k = j \in \{0,1\} \quad \text{ iff } \quad f^k(x)\in N_j.
\]
Then $g$ is a surjection satisfying $g\circ f|_I = \sigma \circ g$ and therefore $f$ is semi-conjugate to the shift map $\sigma$ on $\Sigma_2$.
\end{theorem}

The conjugacy to the model space $\Sigma_2$ may be understood as follows: for any sequence of the symbols $0$ and $1$ there exists an orbit of the discrete system generated by $f$ passing through the sets $N_0$ and $N_1$ in the order given by the sequence. Moreover, if the sequence is $k$-periodic, then so is the orbit.

\begin{corollary}
Let $f$ be a topological horseshoe for $N_0$, $N_1$. Then it follows from Theorem \ref{th:periodic} that for any finite sequence of zeros and ones $(l_0, l_1,\ldots, l_{n-1})$, $l_i\in\{0,1\}$, there exists $x\in N_{l_0}$ such that
\[
f^i(x)\in\inte N_{l_i}
\qquad \text { and } \qquad
f^n(x)=x.
\]
\end{corollary}

\subsection{The theorem on symbolic dynamics for H\'enon 3D map}
In our case, we assume that h-sets $N_i$ are contained in $\mathbb{R}^3$ and the continuous map $f = H^4 : \mathbb{R}^3 \to \mathbb{R}^3$ is such that $N_i$ have exit dimensions equal to 2, that is $u(N_i)=u=2$ and $s(N_i)=s=1$ (as on the right part of Fig. \ref{fig:covering2}).

The sets $a$, $b$ defined in \cite{LiYang} are parallelepipeds spanned by the sets of vertices (Fig. \ref{fig:towel})
\begin{itemize}
\item $a$: \quad $\{(0.84,1.13,0.65), (1.205,1.13,0.65),
(1.205,0.94,1.03), (0.84,0.94,1.03),$ \\
$(0.84,1.01,0.59), (1.205,1.01,0.59),
(1.205,0.82,0.97), (0.84,0.82,0.97) \}$;
\item $b$: \quad $\{(1.365,1.13,0.65),(1.61,1.13,0.65),
(1.61,0.94,1.03), (1.365,0.94,1.03),$ \\
$(1.365,1.01,0.59), (1.61,1.01,0.59),
(1.61,0.82,0.97), (1.365,0.82,0.97) \}$.
\end{itemize}

Define two h-sets $N_a$, $N_b \subset \mathbb{R}^3$ with supports $a$ and $b$, respectively. The supports are images of the cube $B= (N_a)_{c} = (N_b)_{c} =  \overline{\mathbb{B}_2}\times\overline{\mathbb{B}_1} = [-1,1]^3$ by the  following affine transformations:

\begin{align}
|N_a|=a= C_a^{-1}(B)= & \begin{bmatrix}
0.81 \\ 1.0225 \\ 0.975
\end{bmatrix} + \begin{bmatrix}
0.& 0.19& -0.03 \\ 0.1825 & 0.& 0. \\ 0.& -0.095 & -0.06
\end{bmatrix} \cdot B\text{,}
\\
|N_b|=b= C_b^{-1}(B)= & \begin{bmatrix}
0.81 \\ 1.4875 \\ 0.975
\end{bmatrix} + \begin{bmatrix}
0.& 0.19& -0.03 \\ 0.1225 & 0.& 0. \\ 0.& -0.095 & -0.06
\end{bmatrix} \cdot B,
\end{align}

\begin{theorem}\label{th:symbolic}
For h-sets $N_a$, $N_b$ defined above the following chain of covering relations occurs:
\begin{equation}\label{eq:covH4}
	N_a \overset{H^4}{\Longrightarrow} N_a \overset{H^4}{\Longrightarrow} N_b \overset{H^4}{\Longrightarrow} N_b \overset{H^4}{\Longrightarrow} N_a\text{,}
\end{equation}
which proves the existence of symbolic dynamics for the fourth iterate of map $H$, defined by \eqref{eq:h3d}.
\end{theorem}

\begin{proof}

Using the CAPD library for C++ \cite{CAPD} one is able to calculate the (over-estimated) image of the h-set $N$ through the map $f$ and enclose it in a cuboid (an interval closure of $f(N)$, denoted by $[f(N)]$). 

To prove each of four covering relations $N_0 \overset{f}{\Longrightarrow} N_1$, where $f=H^4$, $N_0, N_1 \in \{N_a, N_b\}$ (as in \eqref{eq:covH4}), we check a sufficient condition, which is a conjunction of two:
	\begin{enumerate}[(I)]
		\item\label{enu:1} $[f(N_0)]$ is \emph{spanned across} $N_1$, that is:
		\begin{itemize}
			\item the image projected on any unstable coordinate ($1$ or $2$) lies outside $N_1$: 
		\[
		\left\|[C_{N_1}\circ f\circ C_{N_0}^{-1}(B)]_{1,2}\right\|>1\text{;}
		\]
			\item or else, the image projected on stable coordinate ($3$) lies between the two components of $N_1^+$:
		\[
		\left\|[C_{N_1}\circ f\circ C_{N_0}^{-1}(B)]_3\right\|<1.
		\]
		\end{itemize}
	\item\label{enu:2} The estimated image of the exit set $[f(N_0^-)]$ \emph{lies outside} $N_1$ and it \emph{is homotopy equivalent} to the estimated image of $N_0^-$ through a chosen linear map $(A,0)$.
	
	The linear map $A$ satisfying the conditions from Def. \ref{def:cov} that we choose is $A=D(C_{N_1}\circ f\circ C_{N_0}^{-1})(0)_{u}$. Fig. \ref{fig:lin} compares, as an example, the images of the box $B$ via maps $C_{a}\circ H^4\circ C_{a}^{-1}$ and $(A=D(C_{a}\circ H^4\circ C_{a}^{-1})(0)_{u},0)$. 
	\begin{figure}[h]
	\includegraphics[height=6cm]{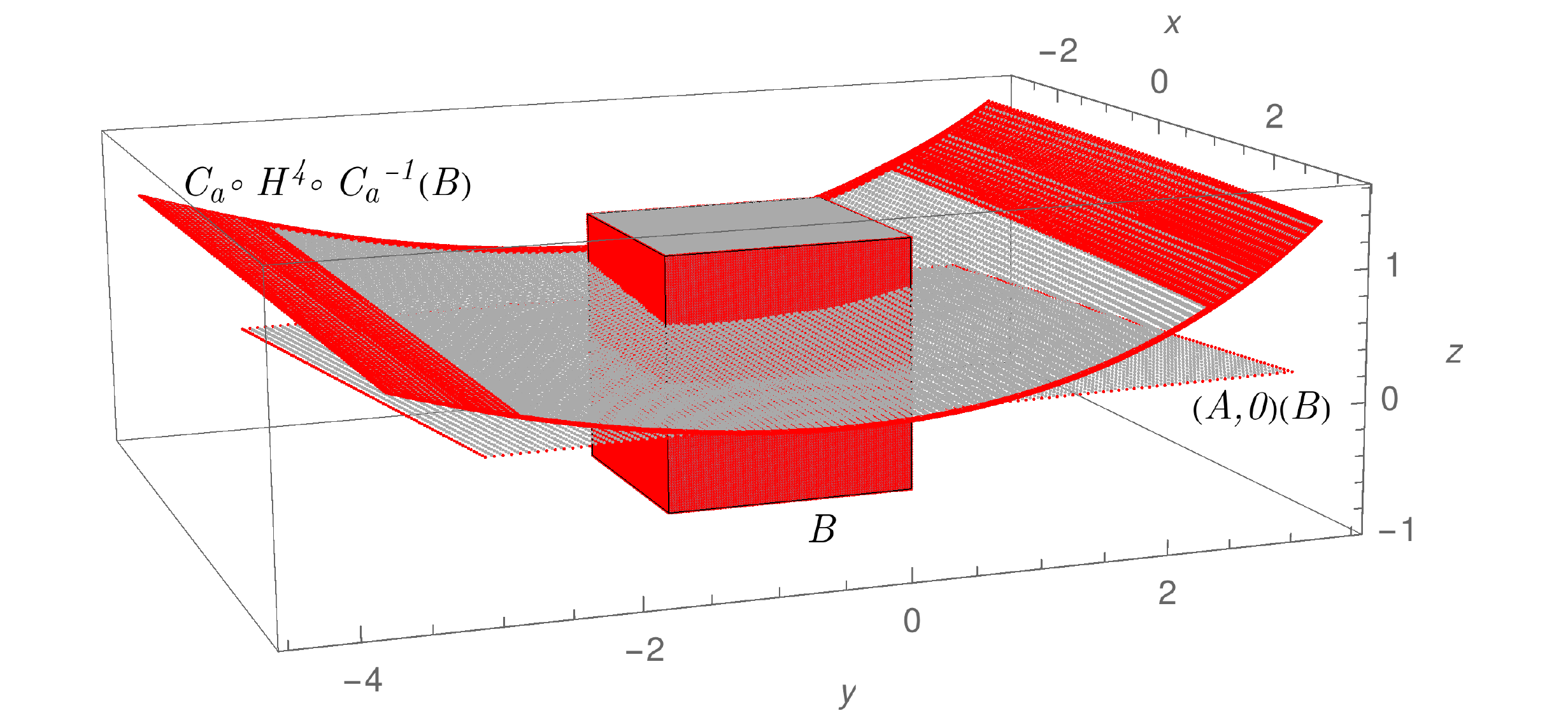}
	\caption{\label{fig:lin}The images of the box $B$ via maps $C_{a}\circ H^4\circ C_{a}^{-1}$ and $(A=D(C_{a}\circ H^4\circ C_{a}^{-1})(0)_{u},0)$. The exit set and its images are marked in red.}
\end{figure}

	Next, we fix one of the unstable dimensions (1 and subsequently 2) and check for each of two faces of $(N_0^-)_{1,2}$, if the interval hull of the union of its image through $f$ and through $(A,0)$ lies outside $N_1$, projected on the considered coordinate:
	\[
		\left\|\big[[C_{N_1}\circ f\circ C_{N_0}^{-1}(B_{\text{(face)}}^-)]\cup [C_{N_1}\circ (A,0)\circ C_{N_0}^{-1}(B_{\text{(face)}}^-)]\big]_{1,2}\right\|>1	
	\]
This means, in particular, that there exists a homotopy
\[
h:[0,1]\times \overline{\mathbb{B}_2}\times\overline{\mathbb{B}_1} \ni (t,x,y) \longmapsto (1-t) f(x,y) + t (A(x),0) \in \mathbb{R}^3
\]
connecting the images of the face through $f$ and through $(A,0)$ which does not touch $N_1$.
\end{enumerate}

\subsubsection*{Outline of the computer-assisted proof \cite{proof}}~

\noindent The program \texttt{03a\_Henon\_Towel\_Symbolic\_Dynamics.cpp} consists of the following steps:

\begin{enumerate}
	\item We define the interval map $H^4$ and h-sets $N_a$, $N_b$.
	\item We check each of four covering relations $N_0 \overset{f}{\Longrightarrow} N_1$, where $f=H^4$, $N_0, N_1 \in \{N_a, N_b\}$ by:
		\begin{enumerate}
			\item dividing $N_0$ in $20\times 20\times 20$ small h-sets and checking the condition \eqref{enu:1} for each part,
			\item dividing each of two faces of $(N_0^-)_1$ in $10\times 10$ small parts and checking the condition \eqref{enu:2} for every part in unstable dimension 1,
			\item dividing each of two faces of $(N_0^-)_2$ in $10\times 10$ small parts and checking the condition \eqref{enu:2} for every part in unstable dimension 2.
		\end{enumerate}
		\item The above conditions are fulfilled, proving every relation $N_0 \overset{H^4}{\Longrightarrow} N_1$, where  $N_0, N_1 \in \{N_a, N_b\}$. Therefore, Theorem \ref{th:symbolic} is proved.
	\end{enumerate}

\end{proof}

\section{Hyperbolicity}

\subsection{Uniform and strong hyperbolicity}
 
First, recall the notion of uniform hyperbolicity.
Let $f:\mathbb{R}^n \to \mathbb{R}^n$ be a diffeomorphism and $I\subset \mathbb{R}^n$ -- a compact invariant set for $f$.

\begin{definition}[\cite{HPS}%, Def. 4
]
We say that $f$ is \emph{uniformly hyperbolic} on $I$ if for every point $x\in I$ the tangent space $T_xI$ is equal to a direct sum $T_xI = E^u_x\oplus E^s_x$ such that
\[
D\!f(x)E^u_x = E^u_{f(x)}\text{,} \qquad D\!f(x)E^s_x = E^s_{f(x)}\text{,}
\]
and for some constants $c>0$, $0<\lambda <1$ independent of $x$ the inequalities
\begin{align*}
\forall\; v\in E^s_x \quad \|D\!f^k(x)v\| < c\lambda^k \|v\|
\text{,}
\\
\forall\; v\in E^u_x \quad \|D\!f^{-k}(x)v\| < c\lambda^k \|v\|
\end{align*}
hold for every $k\geq 0$.
\end{definition}

In \cite{Z} and \cite{W2010} the authors introduce a method for proving hyperbolicity using the notion of \emph{h-set with cones}. In our case, we shall not need this definition except some auxiliary notations.

Denote by $Q$ the $n\times n = (u+s)\times(u+s)$ block matrix
\[
Q = \bmatrix \operatorname{Id}_u & 0 \\ 0 & -\operatorname{Id}_s \endbmatrix
\text{,} 
\]
where by $\operatorname{Id}_k$ we mean the identity matrix of dimension $k\times k$.

Suppose now that we have a set $M=\bigcup_{i=1}^N M_i$, where $M_i\subset \mathbb{R}^n$ are compact and have pairwise disjoint interiors. Each $M_i$ is related to an affine coordinate system $C_i$. Let $I$ be the invariant part of $M$ for the diffeomorphism $f$, that is
\[
I = \operatorname{Inv}_f(M) = \{x\in M \; |\; \forall_{n\in \mathbb{N}}\; f^n(x)\in M \land f^{-n}(x) \in M\}.
\]

Denote also $f_{ij} = C_j\circ f\circ C_i^{-1}$ for $i,j=1,\dots,N$.

\begin{definition}[\cite{W2010}, Def.\ 2.2]

We say that $f$ is \emph{strongly hyperbolic} on $M$ if for $x\in M_i$ and $j=1,\dots,N$ such that $f(M_i)\cap M_j\neq\varnothing$
\begin{equation}\label{cc}
\text{the matrix }\quad D\!f_{ij}(x)^T\cdot Q \cdot D\!f_{ij}(x) - Q \quad
\text{ is positive definite.}
\end{equation}

\end{definition}

\begin{theorem}[\cite{W2010}, Th.\ 2.3]\label{th:hyperbolic}
If $f$ is strongly hyperbolic on $M=\bigcup_{i=1}^N M_i$, then $f$ is uniformly hyperbolic on $I=\operatorname{Inv}_f(M)$.
\end{theorem}

\subsection{The theorem on hyperbolicity}

\begin{theorem}\label{th:myHyperbolic}
The map $H^4$ is uniformly hyperbolic on $\operatorname{Inv}_{H^4}(a\cup b)$.
\end{theorem}

\begin{proof}
	From the computer-assisted proof \cite{proof} we deduce that $H^4$ is strongly hyperbolic on $a\cup b$ and therefore the thesis follows from Theorem \ref{th:hyperbolic}.

\subsubsection*{Outline of the computer-assisted proof  \cite{proof}}~

\noindent The program \texttt{04\_Henon\_Towel\_Hyperbolicity.cpp} consists of the following steps:
\begin{enumerate}
	\item We define four interval maps $f_{aa}=C_a\circ H^4\circ  C^{-1}_a$, and analogously $f_{ab}$, $f_{ba}$, $f_{bb}$. We will check the condition \eqref{cc} for $x\in B = C_a(a)=C_b(b)$.
	\item In our case $u=2$, $s=1$, so we also define
\[
Q = \bmatrix 1 & 0 & 0 \\ 0 & 1 & 0 \\ 0 & 0 & -1 \endbmatrix 
\text{.}
\] 
	\item We calculate the interval upper approximation of the matrices $D\!f(B)^T\cdot Q \cdot D\!f(B) - Q$, $f\in\{f_{aa}, f_{ab}, f_{ba}, f_{bb}\}$, and check with Sylvester's criterion if they are positive definite (the conditions are not fulfilled on the whole box $B$, for the images are too much over-estimated).
	\item We divide the box $B$ into $(25\times 25 \times 25)$ smaller boxes, and for each part $B_i$ and each map $f\in\{f_{aa}, f_{ab}, f_{ba}, f_{bb}\}$ we check if:
	\begin{itemize}
		\item the estimated image $[f(B_i)]$ has an empty intersection with $B$: if so, then $B_i$ does not intersect $I$ and we do not need to check  the condition \eqref{cc};
		\item else, if the interval matrix $[D\!f(B_i)]^T\cdot Q \cdot [D\!f(B_i)] - Q$ is positive definite.
	\end{itemize}
	\item The above conditions are fulfilled for each map $f\in\{f_{aa}, f_{ab}, f_{ba}, f_{bb}\}$ and every point in $B$. Therefore, Theorem \ref{th:myHyperbolic} is proved.
\end{enumerate}

\end{proof}

\section{Conclusion}
We proved the existence of symbolic dynamics in $a\cup b$ for the $\mathbb{R}^3$ map $H^4$ defined by \eqref{eq:h3d}, and also the hyperbolicity of its invariant part $\operatorname{Inv}_{H^4}(a\cup b)$. We believe that the methods used in this paper can be easily applied also in other dimensions, especially for maps with two or more expanding directions.

\section*{Acknowledgment}
The author would like to thank Professors Piotr Zgliczy\'nski and Daniel Wilczak for all their devoted time and
insightful comments during common discussions.

\bibliographystyle{plain} 
\bibliography{hyperchaos_bib}
%\nocite{*}

\end{document}